\documentclass[onefignum,onetabnum]{siamart171218}
\usepackage[margin=4cm]{geometry}
\usepackage{enumitem}

\usepackage{lipsum}
\usepackage{amsfonts}
\usepackage{graphicx}
\usepackage{epstopdf}
\usepackage{algorithmic}
\ifpdf
  \DeclareGraphicsExtensions{.eps,.pdf,.png,.jpg}
\else
  \DeclareGraphicsExtensions{.eps}
\fi

\newsiamremark{remark}{Remark}
\newsiamremark{hypothesis}{Hypothesis}
\crefname{hypothesis}{Hypothesis}{Hypotheses}
\newsiamthm{claim}{Claim}

\headers{Optimal Selection of Transaction Costs}{David Mguni}

\title{Optimal Selection of Transaction Costs in a Dynamic Principal-Agent Problem}

\author{David Mguni\thanks{ 
  (\email{davidmguni@hotmail.com}).}{ Quantitative and Applied Spatial Economic Research Laboratory, University College London, Gower Street, London, WC1E 6BT, UK.\newline Centre for Doctoral Training in Financial Computing \& Analytics, University College London, Gower Street, London, WC1E 6BT, UK. }}
\usepackage{amsopn}

\makeatletter
\newcommand*{\addFileDependency}[1]{
  \typeout{(#1)}
  \@addtofilelist{#1}
  \IfFileExists{#1}{}{\typeout{No file #1.}}
}
\makeatother

\newcommand*{\myexternaldocument}[1]{%
    \externaldocument{#1}%
    \addFileDependency{#1.tex}%
    \addFileDependency{#1.aux}%
}

\theoremstyle{empty}

\newtheorem{refproof}{Proof}[section]

\ifpdf
\hypersetup{
  pdftitle={Optimal Selection of Transaction Costs in a Dynamic Principal-Agent Problem},
  pdfauthor={David Mguni}
}
\fi

\myexternaldocument{ex_supplement}

\begin{document}

\maketitle

\begin{abstract}
  
Environments with fixed adjustment costs such as transaction costs or \lq menu costs\rq$ $ are widespread within economic systems. The presence of fixed minimal adjustment costs produces adjustment stickiness so that agents must choose a sequence of points at which time to perform their actions. This paper performs an analysis of the effect of transaction costs on agent behaviour within a dynamic optimisation problem by way of introducing the theory of incentive design to optimal stochastic impulse control. The setup consists of an agent that maximises their utility by performing a sequence of purchases of some consumable good over some time horizon whilst facing transaction costs and a Principal that chooses the agent's transaction costs. This results in a dynamic Principal-Agent model in which the agent uses impulse controls to perform adjustments to their cash-flow process. We address the question of which fixed value of the transaction cost the Principal must choose to induce a desired behaviour from the agent. We study the effect of changes to the transaction cost and show that with an appropriate choice of transaction cost, the agent's preferences can be sufficiently distorted so that the agent finds it optimal to maximise the Principal's objective even when the agent's cash-flow is unobserved by the Principal.
\end{abstract}



\section{Introduction}
There are numerous environments in which financial agents incur fixed or minimal costs when adjusting their financial positions; trading environments with transaction costs, real options pricing and real estate and large-scale infrastructure investing are a few important examples. However, despite the fundamental relevance in theoretical finance and economic theory, the task of modelling minimally bounded adjustment costs within a dynamic Principal-Agent model, mechanism design or generally strategic interactions with informational asymmetries has as of yet, received no analytic treatment.

In this paper, we analyse the effect of transaction costs within a dynamic Principal-Agent model. In this environment, an agent makes purchases of some costly good over some time horizon. 
Each time the agent performs a purchase, the agent incurs at least some fixed minimal cost (e.g. a transaction cost) which is chosen in advance by a Principal. The cost of each purchase is drawn from the agent's liquidity which is modelled by a jump-diffusion process and is observed only by the agent. When the agent's liquidity process hits $0$, the process is terminated as at this point the agent goes bankrupt. Therefore, since the agent's purchases incur fixed minimal costs, the agent performs a sequence of discrete purchases (possibly of varying size) in order to maximise their utility over the horizon of the problem. Since the Principal gets to choose the fixed value of the transaction cost, the Principal aims to choose a transaction cost that induces a specific consumption behaviour from the agent. Since the agent cannot perform its purchases in a continuous fashion, we model the agent's problem as an impulse control which allows us to study optimal control problems in which each action incurs some fixed cost.  

\subsection*{Overview}
$\indent$The aim of this analysis is twofold: the first objective is to study the effect of introducing a transaction cost on the agent's consumption policy and the relationship between the agent's policy and the transaction cost. The second objective of the paper is to fully determine the value of the transaction cost that induces an agent policy that is desirable for the Principal. Thus in the latter case, the choice of transaction cost serves to condition the agent's preferences so that the timing, magnitude (and total number) of the agent's investment adjustments coincide with the Principal's objectives. The analysis of the paper is performed with sufficient generality to allow for the Principal to be uninformed about the agent's preferences and cash-flow process. Nonetheless, the Principal can transfer wealth to (or from) the agent at the point of the agent's investment adjustments in order to induce desirable changes in the agent's purchasing policy.  \

The analysis of the paper is selected with appeal to investigate financial environments with transaction costs and in which the optimal choice of transaction cost is unknown. The study of  public-private partnerships (e.g. employment initiatives or capital investments within), trading with transaction costs and central authorities that seek to condition the behaviour of players in a given financial environment are some examples. 
\subsection*{Background}
$\indent$Consider firstly the example of a single irreversible investment for a firm that privately observes the demand process. In order to  maximise its overall profit, the firm strategically selects a profit-maximising time to enter the market. Secondly, consider the case of a firm that wishes to adjust its production capacity  according its observations of market (demand) fluctuations in order to maximise its cumulative profits. For the firm, increasing production capacity involves paying investment costs which include fixed costs with which the increases in production yield additional firm revenue. In this case, to maximise overall profits the firm selects some \textit{optimal sequence} of capital adjustments implemented over the firm's time horizon. \

In the case of the  single irreversible firm investment, it is widely known that the optimal firm strategy is to delay investment beyond the point at which the expected returns of investment becomes positive --- from the agent's perspective, the late entry of investment results in a socially inefficient outcome \cite{dixit1994investment}. Similarly, in multiple production capacity case the firm's decision process relating to profit maximising production capital levels often also produce  socially inefficient outcomes. 
\

In both cases, it is therefore natural to ask whether it is possible for an (uninformed) central planner to sufficiently modify the firm's preferences so that the firm's investment decisions produce socially efficient outcomes. The case of a single irreversible firm investment (with asymmetric information) was analysed in \cite{kruse2015optimal} in which it was shown that a regulator can induce socially efficient entry decisions through the use of a posted-price mechanism. \

In particular, in \cite{kruse2015optimal} it is shown that by performing a transfer of wealth at the point of an agent's decision, a central authority or Principal who does not observe the state of the world can sufficiently distort an informed agent's preferences in an optimal stopping problem so that the agent's decision to stop the process coincides with the Principal's optimal stopping time. \

Presently however, the literature concerning multiple sequential investment analysis has been primarily limited to entrance and exit problems within environments of complete information (see for example \cite{Zervos2000}). Thus, the important case of Principal-Agent models with multiple sequential investments has thus far not been studied. \subsubsection*{\textit{Theoretical Framework}}
$\indent$The appropriate modelling framework for multiple sequential investment problem in environments of future uncertainty is optimal stochastic control theory. In stochastic control theory, the inclusion of fixed minimal control costs induces a form of system modifications enacted by the agent or controller known as \textit{impulse control}. Impulse control models are optimal control problems in which the cost of control is bounded below so that modifying the system dynamics incurs at least, some fixed minimum cost. In impulse control models, the dynamics of the system are modified through a sequence of discrete actions or bursts chosen at times that the agent chooses to apply the control policy. This distinguishes impulse control models from the classical (continuous) optimal control models in which players are assumed to continuously make infinitesimally fine adjustments for which the associated costs can be made arbitrarily small.\

Given the discrete nature of the modifying actions of impulse controls, impulse control models represent appropriate modelling frameworks for financial environments with transaction costs, liquidity risks and economic environments in which players face fixed adjustment costs (e.g. \lq menu costs\rq$ $).  More generally, impulse control models are suitable for describing systems in which the dynamics are modified by sequences of discrete, timed actions.\
 
We refer the reader to \cite{bensoussan1982controle} as a general reference to impulse control theory and to \cite{vath2007model, palczewski2010finite} for articles on applications. Additionally, matters relating to the application of impulse control models within finance have been surveyed extensively in \cite{korn1999some}.

\subsubsection*{\textit{Literature}}
$\indent$Current modelling methods of multiplayer  interactions with asymmetric information with multiple ($N>2$) adjustments are modelled by stochastic differential games\footnote{Stochastic differential games represent the multiplayer generalisation of stochastic control theory.} with player controls restricted to those belonging to an  absolutely continuous class of controls (e.g. \cite{cardaliaguet2007differential, cardaliaguet2009stochastic, cardaliaguet2012differential}). In particular, the restriction to absolutely continuous controls implies players modify their positions by performing infinitesimally fine adjustments throughout the horizon of the problem. This renders models with absolutely continuous controls unsuitable for describing behaviour in systems with fixed minimal costs since continuous adjustment would result in immediate ruin. 
\subsection*{Contribution}
The analysis addresses the absence of dynamic Principal-Agent models with fixed minimal costs. Our main result is to determine the value of the transaction cost that induces the Principal's desired consumption policy to be executed by the agent. We also conduct an analysis of the transaction cost parameter and the solution to the agent's optimal control policy.

The results also lead to a solution to the following inverse impulse control problem: 

Let $X^{t_0,x_0}_s=X(s,\omega): [0,T]\times \Omega\to S$ be a one-dimensional diffusion  where $x_0\in S$ and $t_0\in [0,T]$ are parameters that represent the initial point and start time of the process respectively. Suppose that the agent's impulse control problem is specified by the following objective which the agent seeks to maximise by a choice of the control $u\in\mathcal{U}$:
\begin{equation}
J[t_0,x_0;u]=\mathbb{E}\left[\int_{t_0}^{\tau_S}h(s,X_s^{t_{0},x_0,u}) +\sum_{j\geq 1
}c(\tau_{j},z_j)\cdot 1_{\{\tau_j\leq \tau_S\}}+ \phi(X_{\tau_S}^{t_{0},x_0,u})\cdot 1_{\{\tau_S<\infty\}}\right], \label{objective_function_example_principal_agent}
\end{equation}
where $\tau_S:\Omega\to [0,T]$ is some random exit time (i.e. $\tau_S(\omega):=\inf\{s\in [0,T]|X_s^{t_0, x_0  ,\cdot }\in S\backslash A;\;\omega\in\Omega,\; A\subset S\}$ for some measurable subset $A$) and where the control policy takes the form $u(s)=\sum_{j\geq 1}  z_{j}  \cdot 1_{\{ \tau_{j}\leq  T\}}  (s)\in U$ for any $s \in [0,T]$. The quantities $z_1,z_2,\ldots,\in \mathcal{Z}$ and $\tau_1,\tau_2,\ldots$ are $\mathcal{F}-$ measurable intervention times and $\mathcal{F}-$ measurable stopping times where $\mathcal{Z}$ is some admissible set of interventions and $ U$ is a control set. The functions $h:[0,T]\times S\to\mathbb{R}$ and $\phi:S\to\mathbb{R}$ are the running cost and the terminal payoff functions (resp.) where $S\subset\mathbb{R}^q$ is a given fixed domain (solvency region) for some $q\in\mathbb{N}$ and $c: [0,T]\times \mathcal{Z} \to \mathbb{R}$ is an intervention cost function. 

Let  $\mathcal{D}:=\{x\in  S:x<x^{\star}\}$ be a given continuation region, that is, a region in which the agent finds it suboptimal to execute an intervention of any size and suppose there exists an  optimal intervention magnitude $\hat{z}$ that is given by $\hat{z}=\hat{x}-x^{\star}$ for some real-valued constant $\hat{x}$. Lastly, denote by $\lambda\in\mathbb{R}_{>0}$ and $\kappa\in\mathbb{R}_{>0}$ the parameters that represent the \textit{proportional cost} and \textit{fixed cost} parts respectively so that an impulse execution of magnitude $z\in\mathcal{Z}$ incurs a cost $(1+\lambda)z+\kappa$.  The inverse impulse control problem is to determine the value of $\kappa$ and $\lambda$ that induces a given fixed pair $(\hat{x},x^{\star})$ given the objective function $J$ in \eqref{objective_function_example_principal_agent}.

 We perform an analysis of the effect of changes to the parameter $\lambda$ on the quantities $(\hat{x},x^{\star})$. We also determine the values of the fixed cost parameters $\lambda$ and $\kappa$ s.th. given some desired fixed pair of values $(\hat{x},x^{\star})$ $\mathcal{D}\equiv\{x\in  S:x<x^{\star}\}$ and $\hat{z}=\hat{x}-x^{\star}$ that is, we address the question of how to induce a particular impulse control policy through a choice of the transaction costs.

Lastly, our as a corollary to the above theory, we show that the solutions to two distinct optimal impulse control problems can be made to be identical after a transformation that acts purely on the intervention cost function.
\subsection*{Organisation}
$\indent$ In section 1, we give a description of the problem and highlight the connection to optimal stochastic control theory with impulse control. In section 2, we give some definitions central to the apparatus of the impulse control and Principal-Agent problem. In section 3, we give a statement of the main results of the paper which is immediately followed by the main analysis in section 4. We lastly summarise with concluding comments which constitutes Section 5.    

\section{Consumption with Transaction Costs}
Consider an agent that observes its liquidity process (cash-flow) which is subject to exogenous shocks and a Principal that does not observe the process. The agent makes costly purchases and seeks to maximise their consumption over some given time horizon before the point at which the liquidity process hits $0$ (bankruptcy). Each purchase incurs at least some fixed minimal cost or \textit{transaction cost} which is drawn  from the agent's cash-flow. Although the agent observes its own cash-flow, the agent's cash-flow is however not observed by the Principal.

We assume that the market consists of one infinitely divisible good that the agent is able to purchase and consume. The Principal and agent have misaligned payoffs, the Principal however is given the choice of the transaction costs paid by the agent. The Principal therefore aims to choose a fixed value of the transaction cost so as to modify the agent's consumption pattern to satisfy some given objective.

A formal description of the problem is as follows:

\noindent Let $X_s^{t_0,x_0}=X(s,\omega): [0,T]\times \Omega\to S$ be a stochastic process which represents the agent's cash-flow process at a time $s\in [0,T]$ where $t_0\in [0,T]$ and $x_0\in \mathbb{R}$ are parameters that define the start time of the problem and the initial amount of cash held by the agent and $T\in ]0,\infty]$ is the horizon of the problem. When there are no purchases, the agent's cash-flow process evolves according to the following expression:
\begin{align}
X_s^{t_0,x_0}=x_0+\int_{t_0}^{s\wedge \tau_S}\Gamma X_r^{t_0,x_0}dr+\int_{t_0}^{s\wedge \tau_S}\sigma X_r^{t_0,x_0 } dB_r+\int^{s\wedge \tau_S}_{t_0}\int  X_{r-}^{t_0,x_0}\gamma(r,z) \tilde{N}(dr,dz),\nonumber
\\ \mathbb{P}-{\rm a.s},& \nonumber
\\ X^{t_0,x_0}_{t_0}:= x_0,&\label{ch5invliquidityprocess}
\end{align}
where $\tau_S:\Omega\to [0,T]$ is a random exit time or \textit{bankruptcy time} which is defined by $\tau_S(\omega):=\inf\{s\in [0,T]|X_s^{t_0, x_0  ,\cdot }\leq 0\}$ so that $\tau_S$ is the time at which the agent's cash-flow process first hits $0$. The parameter $\Gamma:=r_0+\alpha$ consists of $r_0 \in \mathbb{R}_{>0}$ which is the interest rate and $\alpha \in \mathbb{R}$ which is some constant. The constant $\sigma\in \mathbb{R}$ is the diffusion coefficient and $ S\subset\mathbb{R}$ is the state space. The term $B_r$ is a $1-$dimensional standard Brownian motion and $\tilde{N}(ds,dz)=N(ds,dz)-\nu(dz)ds$ is a compensated Poisson random measure where $N(ds,dz)$ is a jump measure and $\nu(\cdot):= \mathbb{E}[N(1,\cdot)]$ is a L\'{e}vy measure. Both $\tilde{N}$ and $B$ are supported by the filtered probability space and $\mathcal{F}$ is the filtration of the probability space $(\Omega ,\mathbb{P},\mathcal{F}=\{\mathcal{F}_s\}_{s\in [0,T] } )$. We assume that $N$ and $B$ are independent.  

At any time, the agent may make a purchase which incurs some fixed minimal cost. The inclusion of a transaction cost precludes agent control policies for which the agent makes purchases continuously, hence the agent makes purchases over a sequence of times over the horizon of the problem. The sizes of the purchases are $\{z_k \}_{k\in \mathbb{Z}}$ and the sequence of times of the agent's purchases is given by $\{\tau_k (\omega)\}_{k\in \mathbb{N}}$ --- an increasing sequence of $\mathcal{F}_{\tau_k }-$measurable discretionary stopping times so that the agent's control policy is given by the double sequence $(\tau,Z )\equiv \sum_{j\in\mathbb{N}}z_j\cdot1_{\{\tau_j\leq T\}}\in  U$ where $\mathcal{Z}\subset \mathbb{R}$ is the set of feasible agent purchases and $\mathcal{T}$ is a set of  $\mathcal{F}-$measurable stopping times and lastly $U\in\subset \mathcal{T}\times\mathcal{Z}$.

The agent's cash-flow process is therefore affected sequentially at the points of purchases performed by the agent and is described by a stochastic process that obeys the following expression:
\begin{align}\nonumber
X_s^{t_0,x_0,(\tau,Z)}=x_0&+\int_{t_0}^{s\wedge \tau_S}\Gamma X_r^{t_0,x_0,(\tau, Z) }dr-\sum_{j\geq 1} ((1+\lambda)z_j+\kappa)\cdot 1_{\{ \tau_{j}\leq T\}}\\&\begin{aligned}+\int_{t_0}^{s\wedge \tau_S}\sigma X_r^{t_0,x_0,(\tau, Z) } dB_r
+\int^{s\wedge \tau_S}_{t_0}\int  X^{t_0,x_0,(\tau,Z)}_{r-}
\gamma(r,z) \tilde{N}(dr,dz),
\\
 X_{t_0}^{t_0,x_0,\cdot}:=x_0,\;\mathbb{P}-{\rm a.s.}&
\\\label{controlledch2invliquidityprocess}
\forall(t_0,x_0)\in [0,T]\times S,\; \forall s\in [0,T],\;\forall (\tau, Z) \in  U,&
\end{aligned}
\end{align}
where $\kappa, \lambda\in\mathbb{R}_{>0}$ are fixed constants which we shall refer to as the \textit{fixed part} of the transaction cost and \textit{proportional part} of the transaction cost respectively whose pair we denote by $\mathfrak{c}:=(\kappa, \lambda)$. Without loss of generality, we assume that $X^{t_0,x_0,\cdot}_s=x_0$ for any $s\leq t_0$.

The aim of the agent is to maximise their purchases. 

\subsubsection*{Agent Payoff Function}

Given a cash-flow process given by \eqref{controlledch2invliquidityprocess}, the agent's payoff function $ \Pi$ is  given by the following expression:
\begin{align}
\Pi^{(\mathfrak{c},(\tau,Z))} [t_0,x_0]=\mathbb{E}\left[\int_{t_0}^{\tau_S} e^{-\delta r }R(X^{t_0,x_0,(\tau,Z)}_r)dr+\sum_{j\geq 1}  e^{- \delta  \tau_{j} }  c (\tau_{j}^-, z_j ) \cdot 1_{\{ \tau_{j}\leq \tau_S \}}  \right]&,\; \label{franchisEEproffunction_ch_5}
\\\nonumber
\forall (t_0,x_0)\in [0,T]\times S,\forall (\tau,Z)\in & U,
\end{align}
where $R: S\to \mathbb{R}$ is some utility function (we shall later specialise to the case in which $R$ is a power utility function) and $\delta \in ]0,1]$ is the agent's discount factor. The function $c$ is given by $c ( \cdot, z_j )= z_j$ which quantifies the reward endowed to the agent after each purchase. 

In this setting, the Principal chooses a transaction cost which consists of a fixed cost $\kappa \in \mathbb{R}_{>0}$ and a marginal cost parameter $\lambda \in \mathbb{R}_{>0}$ which is proportional to the size of the agent's purchase both of which are incurred by the agent at the point of each purchase. 

The Principal has a payoff function $Q^{(\tau,Z)}$ which is composed of a running gain function $W:[0,T]\times S\to \mathbb{R}$ and a purchase gain function $c_P:[0,T]\times\mathcal{Z}\to\mathbb{R}$. 

\subsubsection*{Principal Payoff Function}

Let $(\tau,Z)\equiv[\tau_j,z_j]_{j\in\mathbb{N}}\in U$ be the agent's policy, then the Principal's payoff function is given by the following:
\begin{align}\label{franchisERproffunction_ch_5}
Q^{(\tau,Z)} [t_0,x_0]=\mathbb{E}\left[\int_{t_0}^{\tau_S}  W(r,X^{t_0,x_0,(\tau,Z)}_r)dr +\sum_{j\geq 1}  e^{- \delta_p  \tau_{j} } c_P (\tau_{j}^-, z_{j} )\cdot 1_{\{ \tau_{j}\leq \tau_S \}}\right]&,
\\\forall (t_0,x_0 )\in [0,T]\times S,\;\forall (\tau,Z)&\in  U,\nonumber
\end{align}
where $W:[0,T]\times S\to\mathbb{R}$ is the Principal's running reward function, the function $c_P:[0,T]\times \mathcal{Z}\to\mathbb{R}$  quantifies the reward endowed to the Principal after each agent purchase and lastly the constant $\delta_p \in ]0,1]$ is the Principal's discount factor. We assume that the Principal purchase gain function $c_P$ is given by $c_P ( \tau_{j}, z_{j} )= \lambda_P  z_{j}+\bar{c}_P  \tau_{j}+ \alpha_P$ where $ \lambda_P,\bar{c}_P, \alpha_P\in \mathbb{R}_{>0} $ are constants.

The agent's problem is to find a sequence of selected magnitudes or an impulse control that alters the agent's cash-flow process in such a way that maximises the agent's payoff.

The problem faced by the Principal is to determine the parameters $(\lambda,\kappa)\in\mathbb{R}_{>0}\times\mathbb{R}_{>0}$ that induce agent purchases at the times and by the magnitudes that the Principal would like (i.e. that coincide with the policy that maximises (\ref{franchisERproffunction_ch_5})), given that the agent seeks to maximise its own objective function \eqref{franchisEEproffunction_ch_5}.

We study the effect of the fixed cost parameters $(\lambda,\kappa)$ associated to the agent's control costs on the agent's consumption pattern. A central aim of this analysis is determine a pair $\mathfrak{c}^\star:=(\kappa^\star,\lambda^\star)$ that maximises the Principal's objective and the conditions under which a desirable agent control policy is induced --- that is, determining the transaction cost that leads to the agent finding it optimal to exercise a control that maximises the Principal's payoff (\ref{franchisERproffunction_ch_5}).

Before embarking on our main analysis, we firstly introduce the stochastic generator associated to the jump-diffusion process. 

The generator of $X$ (of the uncontrolled process) acting on some function $\phi\in \mathcal{C}^{1,2} (\mathbb{R}^l,\mathbb{R}^p)$ is given by: 
\begin{equation} \mathcal{L}\phi(\cdot,x)=\sum_{i=1}^p   \mu_i (x)    \frac{\partial \phi}{\partial x_i}(\cdot,x)+\frac{1}{2} \sum_{i,j=1}^p  (\sigma \sigma^T )_{ij} (x)    \frac{ \partial^2\phi}{\partial x_i\partial x_{j} }+I\phi(\cdot,x),	  \label{generatormain} 
\end{equation}
where $I$ is the integro-differential operator defined by:
\begin{equation} 
I\phi(\cdot,x):= \sum_{j=1}^{l} \int_{\mathbb{R}^p}  \{\phi(\cdot,x+ \gamma^{j}(x,z_j))-\phi(\cdot,x)-\nabla\phi(\cdot,x)  \gamma^{j} (x,z_{j} )\}  \nu_{j} (dz_{j}),\;{\forall  x\in \mathbb{R}^p.} \label{integroopmain}
\end{equation}
\subsubsection*{Controlled State Process} 
The controlled process $X$ which describes the agent's liquidity process is a jump-diffusion process which is affected by impulse controls $u\in  U$. Formally, the agent exercises a control $u(s)=\sum_{j\geq 1}\xi_j \cdot 1_{\{\tau_j\leq T \}}  (s)$ where $0\leq t_0< s\leq T$ and $\xi_1,\xi_2,\ldots\in \mathcal{Z}\subset S$ are impulses that are executed at $\mathcal{F}$-measurable stopping times $\{\tau_i\}_{i\in\mathbb{N}}$ where $0\leq t_0\leq \tau_1< \tau_2< \dots <$ so that an impulse control policy is given by the following double sequence: $u=( \tau_1, \tau_2,\ldots,; z_1, z_2,\ldots,)\in  U$. We assume that the impulses $\xi_j \in \mathcal{Z}$ are $\mathcal{F}-$measurable for all $j \in \mathbb{N}$. Hence, let us suppose that an impulse $\zeta \in \mathcal{Z}$ determined by some admissible policy $u\in U$ is applied at some $\mathcal{F}-$measurable stopping time $\tau:\Omega \to [0,T]$ when the state is $x'=X^{t_0,x_0,\cdot} (\tau^-)$, then the state immediately jumps from $x'=X^{t_0,x_0,\cdot} (\tau^-)$ to $X^{t_0,x_0,u} (\tau)=\Gamma (x',\zeta)$ where $\Gamma :S\times \mathcal{Z}\to S$ is called the impulse response function. 

For any control policy $u=[\tau_j,\xi_j]_{j\in\mathbb{N}}\in U$ and for any $\tau \in \mathcal{T}$, we denote by $ \mu_{[t,\tau]} (u)$ the number of impulses the controller executes within the interval $[t,\tau]$ under $u\in U$ . We say that the impulse control $u\in U$ is admissible on $[0,T]$ if either the number of impulse interventions is finite on average i.e. $\mathbb{E}[\mu_{[0,T]} (u)]<\infty$ or if $\mu_{[0,T]} (u)=\infty\implies\lim_{j\to \infty }\tau_j=\infty$. We shall hereon use the symbol $\mathcal{U}$ to denote the set of admissible impulse controls. 

\subsection*{Notation}
Let $\Omega$  be a bounded open set on  $\mathbb{R}^{p+1}$. Then we denote by:
$\bar{\Omega}$  --- the closure of the set $\Omega$.\\
$Q(s,x;R)={{(s',x' ) \in \mathbb{R} ^{p+1}:\max |s'-s|^{\frac{1}{2}}  ,|x'-x|  }<R,s'<s}$. \\
$\partial \Omega$  --- The parabolic boundary $\Omega$  i.e. the set of points $(s,x) \in \bar{\mathcal{S}}$ s.th. $R>0, Q(s,x;R)\not\subset\bar{\Omega}$.\\
$\mathcal{C}^{\{1,2\}} ([0, T],\Omega )=\{h \in C^{\{1,2\}} (\Omega ): \partial_s h, \partial_{x_i,x_{j} } h \in C(\Omega )\}$, where  $\partial_s$ and  $\partial_{x_{i}, x_{j}}$ denote the temporal differential operator and second spatial differential operator respectively.\\
$\nabla\phi=(\frac{\partial \phi}{\partial x_1 },\ldots,\frac{\partial \phi}{\partial x_p})$ --- The gradient operator acting on some function $\phi \in C^1 ([0,T]\times \mathbb{R}^p)$.\\
$|\cdot|$   --- The Euclidean norm to which $\langle x,y \rangle$  is the associated scalar product acting between two vectors belonging to some finite dimensional space. 

For notational convenience, we use $u=[\tau_j,\xi_j ]_{j\geq 1} $ to denote the agent's control policy $u=\sum_{j\geq 1}\xi_j  \cdot1_{\{\tau_j\leq T \}}  (s)\in \mathcal{U}$. Additionally, where it will not cause confusion and where the time index requires emphasis, we use the notation $X^{t,x_{0}} (s)\equiv X_s^{t,x_{0}}$ for any $ s\in [0,T]$. 

\section{Preliminaries}
\begin{definition}\label{Definition 2.1.}
\noindent The agent and Principal have value functions $v_A$ and $v_P$ that are respectively given by the following expressions:
\begin{equation}
v_A ({t_0,x_0})=\sup_{u\in \mathcal{U}}\Pi^{(\mathfrak{c},u)}[t_0,x_0]  ,\qquad
	v_P ({t_0,x_0})=\sup_{u\in \mathcal{U}}Q^{(u)}[t_0,x_0] , \; \forall (t_0,x_0)\in [0,T]\times S.\label{valuefunctions_prin_agent}
\end{equation}
 \end{definition} 

Where it will not cause confusion, we write $v_A ({t_0,x_0})\equiv v({t_0,x_0})$ for any $(t_0,x_0)\in [0,T]\times S$.

With reference to the Principal's problem (\ref{valuefunctions_prin_agent}), we can express the Principal's problem as the following:

Find $\mathfrak{c}^\star\in\mathbb{R}^2$ s.th.
\begin{equation}
\Pi^{(\mathfrak{c}^\star,u^{\star} )} [{t_0,x_0}]=v_A ({t_0,x_0}),\hspace{7 mm}     Q^{(u^{\star} )} [{t_0,x_0}]=v_P ({t_0,x_0}), \label{implementability2}
\end{equation}

We now give a definition which is central to the problem:

\begin{definition}[Implementability]\label{Definition 2.2. }

We say that $\mathfrak{c}$ implements an impulse control policy\\ $u^{\star}=[ \tau_j^{\star}, z_j^{\star} ]_{j\geq 1}\in \mathcal{U}$ if the following condition is satisfied:
\begin{equation} \label{implementability1_prin_agent}
 \Pi^{(\mathfrak{c},u^{\star} )} [{t_0,x_0}]\geq  \Pi^{(\mathfrak{c},u' )} [{t_0,x_0}], \qquad  {\forall(t_0,x_0)\in [0,T]\times S},\; \forall u'\in \mathcal{U},
\end{equation}
\end{definition}
The implementability condition asserts the optimality of the policy $u^{\star}\in \mathcal{U}$ for the agent, given the transaction cost parameters $\mathfrak{c}$.

Therefore, to analyse the Principal's problem it suffices to characterise $\mathfrak{c}^\star$ and the conditions on the Principal's policy for which the agent always finds it optimal to enact the prefixed impulse control policy $u^{\star}\in  \mathcal{U}$ (so that the inequality in (\ref{implementability1_prin_agent}) is satisfied).

The following object is central to the analysis of impulse control models:
\begin{definition}\label{Definition 2.1.1.}
\

\noindent Let $\tau\in\mathcal{T}$, we define the [non-local] intervention operator $\mathcal{M}:\mathcal{H}\to \mathcal{H}$ acting at a state $X(\tau)$ by the following expression:
\begin{equation}
\mathcal{M} \phi(\tau,X(\tau)):=\inf_{z\in \mathcal{Z}}[\phi(\tau,\Gamma (X(\tau^-),z))+c(\tau, z)\cdot 1_{\{\tau\leq T \}}  ],\label{intervention_operator_definition_equation} \end{equation}
for some function $\phi :[0,T]\times S\to \mathbb{R}$ and $\Gamma : S\times \mathcal{Z}\to  S$ is the impulse response function.
\end{definition}
\section{Main Results}
 
We now present the main results of the paper; we postpone the proofs until the following section. 

\begin{theorem}\label{Theorem 3.1.} 
Let $x^{\star}\in S$ be the Principal's target for the agent's consumption threshold so that whenever the agent's cash flow is less than $x^{\star}$ no purchases are made by the agent. Define $\hat{x}=\hat{z}+x^{\star}$ where $\hat{z}\in\mathcal{Z}$ is the fixed optimal purchase magnitude. Then the agent adopts the Principal's target for the pair $(\hat{x},x^{\star})$ whenever the transaction cost parameter pair $\mathfrak{c}$ is set to the following:
 \begin{align}\nonumber
\lambda^\star(\hat{x},x^{\star})&=  \left(\frac{z}{b}\right)\frac{l_2^{-1}z^{-l_2}-l_1^{-1}z^{-l_1}}{l_1^{-1}z^{-l_1}+l_2^{-1}z^{-l_2}}-1 \\
\kappa^\star(\hat{x},x^{\star})&=z\left[l_1^{-1}+l_2^{-1}-1\right]-z\frac{l_1^{-1}z^{-l_1}-l_2^{-1}z^{-l_2}}{l_1^{-1}z^{-l_1}+l_2^{-1}z^{-l_2}}\left[l_1^{-1}-l_2^{-1}+\ln{\hat{x}}-\ln{x^\star}\right],
 \label{kapparesult33}
\end{align}
where  $b:=\epsilon\delta^{-1}$, $z^m:=\hat{x}^m-x^{\star m}$ and where $\epsilon\in\mathbb{R}/\{0\}$ is a constant that parameterises the agent's risk aversion for the CRRA utility function (c.f. \eqref{principal_agent_crra_utility}) and $\delta$ is the agent's discount factor.

When the agent's liquidity process contains no jumps ($\gamma(z)\equiv 0$ in (\ref{ch5invliquidityprocess})), the parameters $l_1$ and $l_2$ in \eqref{kappa_expression_theorem_principal_agent1} can be expressed exactly in closed form by:
\begin{equation}
l_1=\frac{-1}{b\sigma^2}\left(\sqrt{c^2\delta^2+2b^2\sigma^2\delta}+c\delta\right), \hspace{4 mm} l_2=\frac{1}{b\sigma^2}\left(\sqrt{c^2\delta^2+2b^2\sigma^2\delta}-c\delta\right),
\end{equation}
where $c:=\epsilon\left(\Gamma -\frac{1}{2}\sigma^2\right)$.

For the general case ($\gamma(z)\not\equiv 0$ in (\ref{ch5invliquidityprocess})), the constants $l_1$ and $l_2$ are solutions to the equation:
\begin{equation}
h(l)=0 \label{lequation}
\end{equation}
where the function $h$ is defined by:
\begin{equation}
h(l):=\frac{1}{2}\sigma^2l(l-1)+l\Gamma-\delta+\int_\mathbb{R}\Big\{(1+\gamma(z))^l-1-l\gamma(z)\Big\}\nu(dz).
\end{equation}

If the proportional part $\lambda$ is exogenously fixed, then the value of the fixed part $\kappa$ for which the agent finds it optimal to adopt the Principal's target is given by the following:
 \begin{align}
\kappa^\star(\hat{x},x^{\star},\lambda)=z\left[l_1^{-1}+l_2^{-1}-1\right]+b\left(1+\lambda\right)\left[l_1^{-1}-l_2^{-1}+\ln{\hat{x}}-\ln{x^\star}\right]. \label{kappa_expression_theorem_principal_agent1}
\end{align}
\end{theorem} 
Theorem \ref{Theorem 3.1.} says that if the Principal imposes a transaction cost with proportional part and fixed part given by (\ref{kapparesult33}), then the agent's continuation region is given by $D=\{x< x^{\star}|x,x^{\star}\in S\}$ i.e. the agent makes a purchase whenever the agent's cash-flow attains the value $x^{\star}$. Moreover, the agent's purchase times are $\hat{\tau}_{j+1}=\inf\{s>\tau_j;\; x\geq x^{\star} \}\wedge \tau_S$ and the agent's purchases have a size given by $\hat{z}=\hat{x}-x^{\star}$ which are exactly the intervention times and magnitudes that are optimal for the Principal.

The first result of the theorem relates to the case when the Principal is free to choose the value of the proportional part of the transaction cost parameter $\lambda$ \textit{and} the fixed part of the transaction cost parameter $\kappa$. The second result relates to the case when the Principal is free to choose the value of the fixed part of the transaction cost parameter $\kappa$ but the proportional cost parameter $\lambda$ is exogenous and fixed. 

Theorem \ref{Theorem 3.1.} characterises implementability conditions under which the Principal can sufficiently distort the agent's incentives so that the agent plays actions that maximise the Principal's objective. The following set of results relate to changes in the agent's behaviour following a modification of transaction costs. In particular, the following results characterise the changes in the agent's policy following a change in the agent's transaction costs.

The first result follows from Theorem \ref{Theorem 3.1.}: 

\begin{proposition}\label{Proposition 3.2.}
Let the values $l_1$ and $l_2$ be as in Theorem \ref{Theorem 3.1.} and suppose the initial fixed and proportional costs are given by $\kappa_0\in\mathbb{R}_{>0}$ and $\lambda_0\in\mathbb{R}_{>0}$ respectively. Suppose now that the fixed and proportional costs undergo the transformations $\kappa_0 \to \kappa_1$ and $\lambda_0 \to \lambda_1$, then the agent's intervention threshold and consumption magnitude attain the values $x^{\star}_1=x^{\star}_0+h^\star$ and $\hat{x}_1=\hat{x}_0+\hat{h}$ (respectively) whenever the values $\lambda_1$ and $\kappa_1$ are given by the following expressions:
\begin{align}
\lambda^\star_1(\hat{m},m^\star,\kappa_0, \lambda_0)&=  \left(\frac{\tilde{z}}{b}\right)\frac{l_2^{-1}\tilde{z}^{-l_2}-l_1^{-1}\tilde{z}^{-l_1}}{l_1^{-1}\tilde{z}^{-l_1}+l_2^{-1}\tilde{z}^{-l_2}}-1 \\
\kappa^\star_1(\hat{m},m^\star,\kappa_0, \lambda_0)&=\tilde{z}\left[l_1^{-1}+l_2^{-1}-1\right]\nonumber
\\&-\tilde{z}\frac{l_1^{-1}\tilde{z}^{-l_1}-l_2^{-1}\tilde{z}^{-l_2}}{l_1^{-1}\tilde{z}^{-l_1}+l_2^{-1}\tilde{z}^{-l_2}}\left[l_1^{-1}-l_2^{-1}+\ln{(\hat{m}+\hat{h})}-\ln{(m^\star+h^\star)}\right],
\end{align}
where $\tilde{z}^k:=(\hat{m}+\hat{h})^k-(m^\star+h^\star)^k$ where $\hat{m}$ and $m^\star$ are the solutions to the equations:
\begin{gather}
\mathbf{Q}(\hat{m},m^\star,\kappa_0, \lambda_0)= \left[\begin{array}{c} Q_1(\hat{m},m^\star,\kappa_0,\lambda_0)\\ Q_2(\hat{m},m^\star,\kappa_0,\lambda_0)\\\end{array}\right]=0
\label{Qeqn1}
\end{gather}
where $Q_1$ and $Q_2$ are given by:
\begin{align*}
&\begin{aligned}
Q_1(x,y,q,k):=&\left(l_1x^{l_1}+l_2x^{l_2}\right)\left(y-x-q+b(1+k)[\ln{x}-\ln{y}]\right)
\\&-(x-b(1+k))(y^{l_1}-x^{l_1}+y^{l_2}-x^{l_2}),
\end{aligned}
\\&\begin{aligned}
Q_2(x,y,q,k):=&\left(l_1y^{l_1}+l_2y^{l_2}\right)\left(y-x-q+b(1+k)[\ln{x}-\ln{y}]\right)
\\&-(y-b(1+k))(y^{l_1}-x^{l_1}+y^{l_2}-x^{l_2}).
\end{aligned}
\end{align*}
and where $b:=\epsilon\delta^{-1}$.
\end{proposition}

Proposition \ref{Proposition 3.2.} says that given an initial fixed and proportional cost for the agent, $\kappa_0$ and $\lambda_0$ respectively, a shift of size $h^\star$ and $\hat{h}$ in the agent intervention threshold and consumption magnitudes (respectively) can be induced whenever the fixed and proportional costs are made to be the values $\kappa_1^\star$ and $\lambda_1^\star$ of the proposition. 

Here, interestingly the initial agent intervention threshold $x^\star_0$ and initial consumption magnitude $\hat{x}_0$ do not feature in any of the equations that determine the values $\kappa_1$ and $\lambda_1$, hence the only required data are the shift targets $(h^{\star},\hat{h})$ and the initial cost parameters $(\kappa_0,\lambda_0)$. This is useful for the case in which the Principal does not observe the agent's current consumption threshold and magnitude but seeks to induce a change in those quantities by some given magnitudes.   

Proposition \ref{Proposition 3.2.} tackles instances in which the transaction cost undergoes a transformation. This allows us to compare the agent's behaviour following a switch in transaction cost. The following result analyses the change in the agent's behaviour following (continuous) changes in the transaction costs: 

\begin{proposition}\label{Proposition 3.3.} The marginal rates of change in $\hat{x}(\kappa,\lambda)$ and $x^{\star}(\kappa,\lambda)$ w.r.t. $\lambda$ and $\kappa$ are given by the following expressions:
\begin{align}
\begin{aligned}
&\frac{\partial \hat{x}}{\partial \lambda}&=[f_1(\hat{x},x^{\star})]^{-1},\\
&\frac{\partial x^{\star}}{\partial \lambda}&=[f_2(\hat{x},x^{\star})]^{-1},\\
&\frac{\partial \hat{x}}{\partial \kappa} &=[f_3(\hat{x},x^{\star})]^{-1},\\
&\frac{\partial x^{\star}}{\partial \kappa} &=[f_4(\hat{x},x^{\star})]^{-1}.
\end{aligned}
\end{align}
where the parameters $l_1$ and $l_2$ are solutions to the equation (\ref{lequation}) and the functions $f_1,f_2,f_3$ and $f_4$ are given by (\ref{derivative_1_principal_agent}) - (\ref{derivative_4_principal_agent}).
\end{proposition}

Proposition \ref{Proposition 3.3.} therefore evaluates the change in the intervention threshold and consumption magnitudes due to a marginal change in the cost parameters $\lambda$ and $\kappa$. 

The following corollary follows directly from Theorem \ref{Theorem 3.1.} and relates two general stochastic impulse control problems:

\begin{corollary}\label{Corollary 3.4.} Let $X$ be a stochastic process $X_s=X(s,\omega ): [0,T]\times \Omega\to S$ that evolves according to (\ref{ch5invliquidityprocess}).

Consider the following pair of impulse control problems: 
\renewcommand{\theenumi}{\roman{enumi}}
\begin{enumerate}
\item{ Find $u_1^{\star}=[ \tau_{1_j}^{\star}, z_{1_j}^{\star} ]_{j\in \mathbb{N}}\in \mathcal{U}$ and  $\phi_1\in \mathcal{H}$ s.th.
\begin{equation*} 
 \phi_1 ({t_0,x_0})=J_1^{(u_1^{\star} )} [{t_0,x_0}]=\sup_{u_1\in \mathcal{U}} J_1^{(u_1 )} [{t_0,x_0}],\quad {\forall (t_0,x_0)\in [0,T]\times S}.
\end{equation*} }
\item{Find $u_2^{\star}=[ \tau_{2_j}^{\star}, z_{2_j}^{\star} ]_{j\in \mathbb{N}}\in \mathcal{U}$ and  $\phi_2\in \mathcal{H}$ s.th.
\begin{equation*}
 \phi_2({t_0,x_0})=J_2^{(u_2^{\star})} [{t_0,x_0}]=\sup_{u_2\in \mathcal{U}} J_2^{(u_2 )} [{t_0,x_0}] ,\quad {\forall(t_0,x_0)\in [0,T]\times S}, 
\end{equation*}} 
\end{enumerate}
where the objective functions for problem (i) and (ii) are given by the following expressions:
\begin{align}
 J_1^{(u_1)} [t_0,x_0]&= \mathbb{E}^{[x]}\Bigg[\int_{t_0}^{\tau_S}  \alpha e^{-\delta s}\ln(X_s^{t_0,x_0,u_1 } )ds +\sum_{j\geq 1}  (\lambda_1z_j+\kappa_1)\cdot 1_{\{ \tau_{1_j}\leq \tau_S \}}\nonumber
 \\&\qquad\qquad\qquad\qquad\qquad\qquad\qquad\quad+ \Psi_1 (X_{\tau_S}^{t_0,x_0,u_1 } )\cdot 1_{\{\tau_S <\infty\}}\Bigg],	\label{external_objective_corr_3.4_1}
 \\
J_2^{(u_2 )} [t_0,x_0]&= \mathbb{E}^{[x]}\Bigg[\int_{t_0}^{\tau_S}  F(s,X_s^{t_0,x_0,u_2 })ds +\sum_{j\geq 1}  l_2 (X_{ \tau_{2_{j} }-}^{t_0,x_0,u_2 }, z_{j} )\cdot 1_{\{ \tau_{2_{j} }\leq \tau_S \}} \nonumber
\\&\qquad\qquad\qquad\qquad\qquad\qquad\qquad\quad
+ \Psi_2 (X_{\tau_S}^{t_0,x_0,u_2 })\cdot 1_{\{\tau_S <\infty\}} \Bigg], 	 \label{external_objective_corr_3.4}
\end{align}
where $\alpha \in \mathbb{R}$ and $F,l_2,\Psi_1, \Psi_2$ are bounded Lipschitz continuous functions. Suppose also that the controlled process (with interventions) evolves according to (\ref{controlledch2invliquidityprocess}). Then if $u^{\star}_2\in \arg\hspace{-.3 mm}\sup_{u_2\in \mathcal{U}}j_2^{(u_2 )}(t_0,x_0) ,$ then $u_1^{\star}=u_2^{\star}$ whenever:
\begin{align}
\lambda_1^\star&=  \left(\frac{z_2}{b}\right)\frac{l_2^{-1}z_2^{-l_2}-l_1^{-1}z_2^{-l_1}}{l_1^{-1}z_2^{-l_1}+l_2^{-1}z_2^{-l_2}}-1 \\
\kappa_1^\star&=z\left[l_1^{-1}+l_2^{-1}-1\right]-z_2\frac{l_1^{-1}z_2^{-l_1}-l_2^{-1}z_2^{-l_2}}{l_1^{-1}z_2^{-l_1}+l_2^{-1}z_2^{-l_2}}\left[l_1^{-1}-l_2^{-1}+\ln{\hat{x}_2}-\ln{x^\star_2}\right],
\end{align}
where $\hat{x}_2=x^{\star}_2-z^{\star}_2$ and $z^m_2:=\hat{x}_2^m-x_2^{\star m}$ and the constants $l_1$ and $l_2$ are solutions to the equation $m(l)=0$ where $m$ is defined by:
\begin{equation}
m(l):=\frac{1}{2}\sigma^2l(l-1)+l\Gamma-\delta+\int_\mathbb{R}\Big\{(1+\gamma(z))^l-1-l\gamma(z)\Big\}\nu(dz).\label{1sthequation}
\end{equation}

The parameter $x^{\star}_2\in S$ is the parabolic boundary of the continuation region for problem II, that is to say, given some continuation region for the problem II, $D_2$, each $x^{\star}_2$  is of the form $x^{\star}_2=\{x\in  S:x\in \partial D_2 \}$  and $z^{\star}:= \arg\hspace{-0.55 mm}\sup_{z\in \mathcal{Z}}\{\phi_2(\tau_k,\Gamma (X(\tau_k-),z))+l_2(X(\tau_k),z)\}$ quantifies the optimal intervention magnitude for the problem with payoff function $J_2$.
\end{corollary}

Corollary \ref{Corollary 3.4.} says that impulse control problem I has the same optimal control policy solution as that of problem II whenever the intervention cost function in \eqref{external_objective_corr_3.4_1} has a proportional cost and fixed cost given by $\lambda_1^\star$ and $\kappa_1^\star$ respectively. 

\section{Main Analysis}
 We begin by proving Theorem \ref{Theorem 3.1.} which is demonstrated by showing that given $\mathfrak{c}^\star:=(\kappa^\star,\lambda^\star)$ defined in (\ref{kapparesult33}), it is optimal for the agent to execute the sequence of  interventions that maximises the Principal's payoff $Q$.

Before deriving the main results, we require some background results. In particular, we require a verification theorem for the single controller optimal stochastic control problem which was reported in Theorem 6.2 in \cite{oksendalapplied2007}:
\begin{theorem}[Theorem 6.2 in \cite{oksendalapplied2007}]\label{corollary_verification_theorem_zero_sum_degenerate}

Consider the impulse control problem in which the dynamics under the influence of impulse controls $u=[\tau_j,\xi_j]_{j\geq 1}\in\mathcal{U}$ evolves according to the jump-diffusion process $\forall r\in [0,T];\;\forall (t_0,x_0)\in [0,T]\times S$:
\begin{align*}\nonumber
\hspace{-1.8 mm}X_r^{t_0,x_0,u}=x_0+\int_{t_0}^{r}\mu(s,X^{t_0,x_0,u}_s)ds+\int_{t_0}^{r}\sigma(s,X^{t_0,x_0,u}_s)dB_s
+\sum_{j\geq 1}\xi_j  \cdot 1_{\{\tau_j\leq r\}}  (r)&
\\+\int_{t_0}^{r}\int\gamma (X^{t_0,x_0,u}_{s-},z) \tilde{N}(ds,dz),& 
\\\mathbb{P}-{\rm a.s.},&
\end{align*}
and for which the agent seeks to maximise the following objective function:
\begin{align}\nonumber 
&\hspace{-9 mm}J [t_0, x_0  ;u]
\\&\begin{aligned}= \mathbb{E}\left[\int_{t_0}^{\tau_s} f (s, X_s^{t_0, x_0  ,u } ) ds  + \sum_{m\geq 1}  c (\tau_m  , \xi _m  )  \cdot 1_{\{\tau_m  \leq  \tau_S \}}+G (\tau_S , X_{\tau_S  }^{t_0, x_0  ,u } )1_{\{\tau_S  <\infty\}}\right],&\\ \forall (t_0,x_0)\in[0,T]\times S.&
\end{aligned}
\end{align}

Suppose that there exists a function $\phi\in \mathcal{C}^{1,2} ([0,T],S)\cap\mathcal{C}([0,T],\bar{S})$ that satisfies technical conditions (T1) - (T4) and the following conditions:
\renewcommand{\theenumi}{\Roman{enumi}}
 \begin{enumerate}[leftmargin= 6 mm]
	\item  $\phi\leq \mathcal{M}\phi$ on $S$ and define the region  $D$ by:\\
$D=\{x \in S;\phi(\cdot,x)<\mathcal{M}\phi(\cdot,x)\}$ 
where $D$ is the controller continuation region and where $\mathcal{M}$ is the (non-local) intervention operator defined in \eqref{intervention_operator_definition_equation}.

	\item $ \frac{\partial \phi}{\partial s}+\mathcal{L}\phi(\cdot,X^{\cdot,u } (\cdot))+f(\cdot,X^{\cdot,u } (\cdot))\geq 0,\hspace{1 mm}   \forall  u \in \mathcal{U}$ on $S\backslash{\partial D}$.
	
	\item \label{hjb_equation_single_ipulse_corr} $\frac{\partial \phi}{\partial s}+\mathcal{L}\phi(\cdot,X^{\cdot,\hat{u}} (\cdot))+f(\cdot,X^{\cdot,\hat{u}} (\cdot))=0 \hspace{1 mm}$ in $D. $
	\item  $X^{\cdot,u} (\tau_S ) \in \partial S$, $ \mathbb{P}-{\rm a.s}$. on ${\tau_S< \infty}$ and $\phi(s,X^{\cdot,u} (s))\to G(\tau_S,X^{\cdot,u} (\tau_S))$  as $s\to \tau_S^-  $ $\mathbb{P}-{\rm a.s}.,\forall u \in \mathcal{U}$.
\end{enumerate}
Put $\hat{\tau}_0\equiv t_0$ and define $\hat{u}:=[\hat{\tau}_j,\hat{\xi}_j ]_{j \in \mathbb{N}}$ inductively by:\\ $\hat{\tau}_{j+1}= \inf\{s>\tau_j;X^{\cdot,\hat{u}_{[t_0,s]}} (s)\notin D \}\wedge\tau_S$,
then $ \hat{u}\in\mathcal{U}$ is an optimal control for the agent's impulse control problem, that is to say we have:\begin{equation}
\phi(t_0,x_0)= \inf_{u \in \mathcal{U}}J [t_0, x_0  ;u]=J [t_0, x_0  ;\hat{u}];	
\qquad{\forall (t_0,x_0) \in [0,T]\times S}.
\end{equation}
\end{theorem}

\begin{remark}
\label{Remark 4.3}
Let us denote by $\mathcal{D}$ the region $\mathcal{D}=\{x \in S:v (\cdot,x)< \mathcal{M}v (\cdot,x)\}$ so that $\mathcal{D}$ represents the region in which the agent finds an immediate intervention suboptimal.  We can infer the existence of a value $x^{\star}\in S$ for which $\partial \mathcal{D}=\{X_s^{\cdot}=x^{\star}|x^{\star}\in S, s\in[0,T]\}$, that is to say the agent performs an intervention as soon as the cash-flow process $X$ attains a value $x^{\star}$, hence we shall hereon refer to the value $x^{\star}$ as the agent's \text{intervention threshold}.
\end{remark}

\begin{refproof}[Proof of Theorem \ref{Theorem 3.1.}] 

We now seek to characterise the cost function parameters $\lambda$ and $\kappa$ which implement the Principal's control policy. 

Suppose that the agent makes purchases according to the policy $[\tau_k,z_k]_{k\geq 1}\equiv(\tau,Z)\in \mathcal{U}$, hence the agent's payoff function is given by the expression:
\begin{equation*}
\Pi^{(\mathfrak{c},(\tau,Z))} [t_0,x_0]=\mathbb{E}\left[\int_{t_0}^{\tau_S}e^{-\delta r }R(X^{t_0,x_0,(\tau,Z)}_r)dr+\sum_{j\geq 1}  e^{- \delta  \tau_{j} } z_j\cdot 1_{\{ \tau_{j}\leq T \}}  \right], 
\end{equation*}
Let us define the control $(\tau^{\star},Z^{\star})\in\mathcal{U}$ by the following construction: 
\begin{equation*}
\Pi^{(\mathfrak{c},(\tau^{\star},Z^{\star}))} [s,x ]=\sup_{(\tau,Z)\in\mathcal{U}}\Pi^{(\mathfrak{c},(\tau,Z))} [s,x], \; \forall(s,x)\in[0,T]\times S,
\end{equation*}
so that given some $\mathfrak{c}\in\mathbb{R}^2$, the agent's optimal purchase strategy is given by $(\tau^{\star},Z^{\star})\in\mathcal{U}$.\

Recall that the state process obeys the following:
\begin{align}
X_s^{t_0,x_0,(\tau,Z)}=x_0&+\int_{t_0}^{s\wedge \tau_S}\Gamma X_r^{t_0,x_0,(\tau, Z) }dr-\sum_{j\geq 1} ((1+\lambda)z_j+\kappa)\cdot 1_{\{ \tau_{j}\leq \tau_S \}}\nonumber
\\&\begin{aligned}
+\int_{t_0}^{s\wedge \rho}\sigma X_r^{t_0,x_0,(\tau, Z) } dB_r\label{stateprocess2}+\int^s_{t_0}\int  X^{t_0,x_0,(\tau,Z)}_{r-}\gamma(r,z) \tilde{N}(dr,dz),& \\\mathbb{P}-{\rm a.s.}\;  X_{t_0}^{t_0,x_0}:= x_0,&
\\
\forall(s,x),(t_0,x_0)\in [0,T]\times S,\forall (\tau, Z) \in \mathcal{U}.&
\end{aligned}
\end{align}
We now specialise to the case in which the agent's utility function $R$ is given by:
\begin{equation}
R(x)=\epsilon\ln{x},\label{principal_agent_crra_utility}
\end{equation} for some constant $\epsilon \in \mathbb{R}\backslash\{0\}$ so that $R$ can be viewed as a limiting case of the CRRA utility function that is $R(x)=\underset{\eta\to 1}{\lim}\epsilon\frac{x^{1-\eta}-1}{1-\eta} $. 


We note also that given some test function $\phi\in\mathcal{C}^{\{1,2\}}([0,T],\mathbb{R})$, the generator $\mathcal{L}$ for (\ref{stateprocess2}) is given by the following expression (c.f. (\ref{generatormain})): 
\begin{align}
\mathcal{L}\phi(s,x)=\Gamma x\frac{\partial \phi}{\partial x}(s,x) 
&+\frac{1}{2}\sigma^2x^2\frac{\partial^2 \phi}{\partial x^2}(s,x) 
\\&\begin{aligned}+ \int_\mathbb{R}\Big\{\phi(s,x(1+\gamma(z))-\phi(s,x)-x\gamma(z)\frac{\partial\phi}{\partial x}\Big\}\nu(dz),\label{generatorspecific}&
\\ \forall (s,x)\in\mathbb{R}_{>0}\times\mathbb{R}.&\nonumber
\end{aligned}
\end{align}
By (\ref{hjb_equation_single_ipulse_corr}) of Theorem \ref{corollary_verification_theorem_zero_sum_degenerate}, we have that on ${D}$ the following expression holds:
\begin{equation}
R+\frac{\partial \phi}{\partial s} +\mathcal{L}\phi =0. \label{continutityequation}
\end{equation}
Hence, using (\ref{generatorspecific}) and by (\ref{continutityequation}) we have that:
\begin{align}\nonumber
0=e^{-\delta s}\epsilon\ln{x}+\frac{\partial \phi}{\partial s}(s,x)&+\Gamma x\frac{\partial \phi}{\partial x}(s,x) +\frac{1}{2}\sigma^2x^2\frac{\partial^2 \phi}{\partial x^2}(s,x) 
\\&+ \int_\mathbb{R}\Big\{\phi(s,x(1+\gamma(z))-\phi(s,x)-x\gamma(z)\frac{\partial\phi}{\partial x}\Big\}\nu(dz). \label{contequationexplicit1}
\end{align}
Let us try the following ansatz for the candidate function for $\phi$:
\begin{equation}
\phi\equiv\phi_a+\phi_b, \label{phiansatz}
\end{equation} where
\begin{equation}
\phi_a(s,x)=e^{-\delta s}ax^l, \hspace{3 mm} \phi_b(s,x) =e^{-\delta s}\left(b\ln{x} +c\right).    
\end{equation}
for some constants $a,b,c \in\mathbb{R}$.

We firstly seek to ascertain the values of the constants $a,b$ and $c$ hence, inserting the expression for $\phi$ into (\ref{contequationexplicit1}) we find that:
\begin{equation}
h_a(l)+h_b(x)=0,
\end{equation}
where the functions $h_a$ and $h_b$ are given by:
\begin{gather}h_a(l)=\frac{1}{2}\sigma^2l(l-1)+l\Gamma-\delta+\int_\mathbb{R}\Big\{(1+\gamma(z))^l-1-l\gamma(z)\Big\}\nu(dz)\label{contequationexplicit1.2}\\h_b(x)=
\epsilon \ln{x}-\delta(b\ln{x}+c)+b\Gamma-\frac{1}{2}\sigma^2b+\int_\mathbb{R}\big\{b\ln(1+\gamma(z))-b\gamma(z)\big\}\nu(dz), \label{contequationexplicit2}
\end{gather}
from which we find that the equation $h_b(x)=0$ is solved by the following values for $b$ and $c$:
\begin{gather}
b=\epsilon\delta^{-1},\label{b}\\ c=\epsilon\delta^{-2}\left(\Gamma -\frac{1}{2}\sigma^2\right)+\epsilon\delta^{-2}\int_\mathbb{R}\big\{\ln(1+\gamma(z))-\gamma(z)\big\}\nu(dz).
\label{c}
\end{gather}
Let us make a brief excursion to discuss the case when the process (\ref{ch5invliquidityprocess}) contains no jumps i.e. when $\gamma\equiv 0$. In this case, we readily observe that the constants $b$ and $c$ are given by:
\begin{align}
b&=\epsilon\delta^{-1},\\ c&=\epsilon\delta^{-2}\left(\Gamma -\frac{1}{2}\sigma^2\right).
\end{align}

Additionally, (\ref{contequationexplicit1.2}) now reduces to the following expression:
\begin{equation}
h_{a,0}(l):=h_{a}(l)\big|^{\gamma\equiv 0}=\frac{1}{2}\sigma^2l^2+c\delta b^{-1}l-\delta.
\end{equation}
After some simple algebra, we then deduce that in this case there exist two solutions to the equation $h_{a,0}(l)=0$, namely $l_{1,0}$ and $l_{2,0}$ given by:
\begin{equation}
l_{1,0}=\frac{-1}{b\sigma^2}\left(\sqrt{c^2\delta^2+2b^2\sigma^2\delta}+c\delta\right) , \hspace{4 mm} l_{2,0}=\frac{1}{b\sigma^2}\left(\sqrt{c^2\delta^2+2b^2\sigma^2\delta}-c\delta\right). 
\end{equation}
 Let us now return to the case when the process (\ref{ch5invliquidityprocess}) contains jumps. Using (\ref{contequationexplicit1.2}), we now make the following observations:
\begin{equation}
\lim_{m\to\infty}h_a(m)=+\infty, \hspace{5 mm} h_a(m)\big\vert^{m=0}=-\delta.
\end{equation}
Hence, we deduce the existence of values $l_1,l_2$ s.th.
\begin{equation}
h_a(l_1)=h_a(l_2)=0.
\end{equation}
W.l.o.g. let us assume that $l_1<l_2$, since $\forall$ $l,z$ we have that: $(1+\gamma(z))^l-1-l\gamma(z)\}\nu(dz)>0$ so that:
\begin{equation}
|l_1|>l_1,
\end{equation}
and
\begin{equation}
l_1<0<l_2.
\end{equation}
Therefore, the function $\phi$ is given by the following (c.f. (\ref{phiansatz})):
\begin{equation}
\phi(s,x)=e^{-\delta s}[a_1x^{l_1}+a_2x^{l_2}+b\ln{x}+c],\label{phieqnfinal}
\end{equation} where $a_1$ and $a_2$ are a pair of as of yet, undetermined constants and $b$ and $c$ are given by (\ref{b}) - (\ref{c}) and $l_1$ and $l_2$ are solutions to \eqref{contequationexplicit1.2}.

 Our ansatz for the continuation region ${D}$ is that it takes the form:
\begin{equation}
D=\{x< x^{\star}|x,x^{\star}\in S\}. 
\end{equation}
We now seek to determine the value of $x^{\star}$ and characterise the optimal intervention magnitude $\hat{z}$.\\
Now by Theorem \ref{corollary_verification_theorem_zero_sum_degenerate}, we find that for all $x_1\geq x^{\star}$ we have:
\begin{equation}
\phi(\cdot,x)=\mathcal{M}\phi(\cdot,x)=\sup_{z\in\mathcal{Z}}\{\phi(\cdot,x-\kappa-(1+\lambda)z)+z)\}.\label{interventioneqduop_ch_5}
\end{equation}
    We wish to determine the value $z$ that maximises (\ref{interventioneqduop_ch_5}), hence let us now define the function $G$ by the following expression:
\begin{equation}
G(t,z)=\phi(t,x-\kappa-(1+\lambda)z)+z,\quad \forall\; z\in\mathcal{Z}, \forall (t,x) \in [0,T]\times S. \label{functionG_ch_5}\end{equation}
Our task now is to evaluate the maxima of (\ref{functionG_ch_5}) from which we readily observe that the first order condition for the maximum of $G$ is given by:
\begin{equation}
\phi'(\cdot,x-\kappa-(1+\lambda)\hat{z})=\frac{1}{1+\lambda}. 
\end{equation}
Let us now consider a unique point $\hat{x} \in ]0,x^{\star}[$ then using (\ref{interventioneqduop_ch_5}) we find that:
\begin{equation}
\phi'(\cdot,\hat{x})=\frac{1}{1+\lambda}. 
\end{equation}
We now observe that the following expression holds:
\begin{equation}
x^{\star}-\kappa-(1+\lambda)\hat{z}=\hat{x}. 
\end{equation}
We now find that:
\begin{equation}
\hat{z}(x)=\frac{x-\hat{x}-\kappa}{(1+\lambda)}.
\end{equation}
We therefore deduce that $\phi$ is given by the following expression $\forall x \in S$:
\begin{equation}
\phi(\cdot,x)=\phi(\cdot,\hat{x})+\hat{z}. \label{interventioneqdiff}
\end{equation}
 Using (\ref{interventioneqduop_ch_5}) - (\ref{interventioneqdiff}) we readily obtain the following equations:
\begin{align}
\phi'(\cdot,\hat{x})&=\frac{1}{1+\lambda}\label{psiderivativeeqn1.1}\\
\phi'(\cdot,x^{\star})&=\frac{1}{1+\lambda}
\\
\phi(\cdot,x^{\star})-\phi(\cdot,\hat{x})&=\frac{x^{\star}-\hat{x}-\kappa}{1+\lambda}.\label{psieqn_ch_5}
\end{align}
We now separate the analysis into two cases; case I in which the proportional part of the transaction cost $\lambda$ is fixed and, case II in which the Principal is free to choose both values $\kappa,\lambda$.

\noindent\textbf{Case I}

\noindent Inserting
(\ref{phieqnfinal}) into (\ref{psiderivativeeqn1.1}) - (\ref{psieqn_ch_5}) and by the high contact principle\footnote{Recall that the high contact principle is a condition that asserts the continuity of the value function at the boundary of the continuation region.}, we arrive at the following system of equations:
 \renewcommand{\theenumi}{(\roman{enumi})}\begin{enumerate}[leftmargin= 4 mm]
\item\hspace{40 mm} $a_1l_1\hat{x}^{l_1-1}+a_2l_2\hat{x}^{l_2-1}+\frac{b}{\hat{x}}=\frac{1}{1+\lambda} \label{impulseitem1}$
\item \hspace{40 mm} $a_1l_1x^{*l_1-1}+a_2l_2x^{*l_2-1}+\frac{b}{x^{\star}}=\frac{1}{1+\lambda}\label{impulseitem2}$
\item \hspace{40 mm} $a_1(x^{*l_1}-\hat{x}^{l_1})+a_2(x^{*l_2}-\hat{x}^{l_2 })=\frac{x^{\star}-\hat{x}-\kappa}{1+\lambda}+b\ln\left(\frac{\hat{x}}{x^{\star}}\right)\label{impulseitem3}$
\end{enumerate} where $b:=\epsilon\delta^{-1}$.

\noindent The system of 3 equations \ref{impulseitem1} - \ref{impulseitem3} contains 3 unknowns ($\kappa, a_1, a_2$), hence we can solve it, in particular using \ref{impulseitem1} - \ref{impulseitem2} to solve for $a_1$ and $a_2$ we find that:
\begin{gather}
a_1=\left[\frac{z}{1+\lambda}+b\right]l_1^{-1}z^{-l_1}, \label{a_1}\\ a_2=\left[\frac{z}{1+\lambda}-b\right]l_2^{-1}z^{-l_2} \label{a_2},
\end{gather}where $b:=\epsilon\delta^{-1}$ and $z^m:=\hat{x}^m-x^{\star m}$.

After substituting (\ref{a_1}) - (\ref{a_2}) into \ref{impulseitem3} we readily obtain the expression for the fixed cost parameter $\kappa$:
 \begin{align}
\kappa(\hat{x},x^{\star},\lambda)=z\left[l_1^{-1}+l_2^{-1}-1\right]+b\left(1+\lambda\right)\left[l_1^{-1}-l_2^{-1}+\ln{\hat{x}}-\ln{x^\star}\right],
\end{align}
 Hence, given a pair of target cost parameters $({x}^{\star},\hat{{x}})$, we see that any optimal control for the agent's intervention threshold becomes $\tilde{x}^{\star}$ and optimal consumption magnitude becomes $z=x^{\star}-\hat{{x}}$.

\noindent\textbf{Case II}

\noindent We now seek to identify the parameters $\kappa$ and $\lambda$, after setting $a_1=a_2:=a$ in (\ref{phieqnfinal}), then substituting  into (\ref{psiderivativeeqn1.1}) - (\ref{psieqn_ch_5}) we arrive at the following system of equations:
 \renewcommand{\theenumi}{(\roman{enumi})}\begin{enumerate}[leftmargin= 4 mm]
\item\hspace{40 mm} $a(l_1\hat{x}^{l_1-1}+l_2\hat{x}^{l_2-1})+\frac{b}{\hat{x}}=\frac{1}{1+\lambda} \label{impulseitem1.1}$
\item \hspace{40 mm} $a(l_1x^{*l_1-1}+l_2x^{*l_2-1})+\frac{b}{x^{\star}}=\frac{1}{1+\lambda}\label{impulseitem2.2}$
\item \hspace{40 mm} $a(x^{*l_1}-\hat{x}^{l_1}+x^{*l_2}-\hat{x}^{l_2 })=\frac{x^{\star}-\hat{x}-\kappa}{1+\lambda}+b\ln\left(\frac{\hat{x}}{x^{\star}}\right)\label{impulseitem3.3}$
\end{enumerate} where the constant $b$ is given by Equation (\ref{b}).\

The system \ref{impulseitem1.1} - \ref{impulseitem3.3} which involves 3 equations now consists of 3 unknowns ($\kappa, \lambda, a$), hence we can solve for the three unknown parameters.  Eliminating the constant $a$ from the system \ref{impulseitem1.1} - \ref{impulseitem3.3} yields the following expressions for the cost parameters:
 \begin{align}
\lambda(\hat{x},x^{\star})&=  \left(\frac{z}{b}\right)\frac{l_2^{-1}z^{-l_2}-l_1^{-1}z^{-l_1}}{l_1^{-1}z^{-l_1}+l_2^{-1}z^{-l_2}}-1 \label{lambda_case_ii_principal_agent} \\
\kappa(\hat{x},x^{\star})&=z\left[l_1^{-1}+l_2^{-1}-1\right]-z\frac{l_1^{-1}z^{-l_1}-l_2^{-1}z^{-l_2}}{l_1^{-1}z^{-l_1}+l_2^{-1}z^{-l_2}}\left[l_1^{-1}-l_2^{-1}+\ln{\hat{x}}-\ln{x^\star}\right],\label{kappa_case_ii_principal_agent}
\end{align}
where $z^m:=\hat{x}^m-x^{\star m}$.

Hence, given a pair of target cost parameters $({x}^{\star},\hat{{x}})$, we see that any optimal control for the agent's intervention threshold becomes $\tilde{x}^{\star}$ and optimal consumption magnitude becomes $z=x^{\star}-\hat{{x}}$.\

By inverting the procedure, we can further deduce that using equations \ref{impulseitem1.1} - \ref{impulseitem3.3}, we can derive the values $m^\star$, $\hat{m}$ s.th.
\begin{align}
x^{\star}&=m^\star(\kappa,\lambda) \label{f1}\\
\hat{x}&=\hat{m}(\kappa,\lambda), \label{f2}
\end{align}
where $\hat{m}$ and $m^\star$ are solutions to the system of equations:
\begin{equation}
\mathbf{Q}(\hat{m},m^\star,\kappa, \lambda)= \left[\begin{array}{c} Q_1(\hat{m},m^\star,\kappa,\lambda)\\ Q_2(\hat{m},m^\star,\kappa,\lambda)\\\end{array}\right]=0
\end{equation}
where $Q_1$ and $Q_2$ are given by:
\begin{align}
Q_1(x,y,q,k)&:=\left(l_1x^{l_1}+l_2x^{l_2}\right)\left(y-x-q+b(1+k)[\ln{x}-\ln{y}]\right)\nonumber
\\&\qquad-(x-b(1+k))(y^{l_1}-x^{l_1}+y^{l_2}-x^{l_2}),\label{q_1_principal_agent}\\
Q_2(x,y,q,k)&:=\left(l_1y^{l_1}+l_2y^{l_2}\right)\left(y-x-q+b(1+k)[\ln{x}-\ln{y}]\right)\nonumber
\\&\qquad-(y-b(1+k))(y^{l_1}-x^{l_1}+y^{l_2}-x^{l_2}).\label{q_2_principal_agent}
\end{align}
where $l_1$ and $l_2$ are solutions to \eqref{contequationexplicit1.2} and $b:=\epsilon\delta^{-1}$.
$\hfill \square$
\end{refproof}
Though it is not possible to obtain a closed analytic solution to (\ref{q_1_principal_agent}) - (\ref{q_2_principal_agent}), the values $\hat{m}$ and $m^\star$ can be approximated using numerical methods.

Having proven Theorem \ref{Theorem 3.1.}, we can straightforwardly prove Proposition \ref{Proposition 3.2.}:
\begin{refproof}[Proof of Proposition \ref{Proposition 3.2.}]
To prove Proposition \ref{Proposition 3.2.}, we firstly note that using \eqref{f1} - \eqref{f2}, we can express the unobservable parameter pair $(x^\star_0,\hat{x}_0)$ in terms of the observable parameter pair $(\tilde{\lambda}_0,\kappa_0)$, that is $x^{\star}_0=m^\star(\kappa_0,\lambda_0)$ and $\hat{x}_0=\hat{m}(\kappa_0,\lambda_0)$. Hence, we have that $x^{\star}_1=m^\star(\kappa_0,\lambda_0)+h^\star$ and $\hat{x}_1=\hat{m}(\kappa_0,\lambda_0)+\hat{h}$, where $x^\star_1$ and $\hat{x}$ are the target consumption level and target consumption threshold respectively. Inserting these expressions for $x^{\star}_1$ and $\hat{x}_1$ into \eqref{lambda_case_ii_principal_agent} and \eqref{lambda_case_ii_principal_agent} yields the result.
$\hfill \square$
\end{refproof}
We have therefore succeeded in providing a full characterisation of the parameters of transaction costs that sufficiently distort the incentives of a rational agent so that the agent finds it optimal to maximise the Principal's payoff. In particular, if the above values for the transaction cost are adopted by the Principal, the rational agent finds it optimal to adopt a consumption pattern that is optimal for the Principal. 

We now give a sketch of the remaining proofs, the first of which follows from direct calculation:
\begin{refproof}[Proof of Proposition \ref{Proposition 3.3.}]

\noindent To prove Proposition \ref{Proposition 3.3.}, we differentiate \ref{impulseitem1.1} and \ref{impulseitem2.2} w.r.t. $\hat{x}$ and $x^{\star}$ respectively and plugging in (\ref{f1}) and (\ref{f2}).\\ We now observe that $\frac{\partial \hat{x}}{\partial \lambda},\frac{\partial x^{\star}}{\partial \lambda},\frac{\partial \hat{x}}{\partial \kappa},\frac{\partial x^{\star}}{\partial \kappa}$ are given by the following expressions:
\begin{align}
\frac{\partial \hat{x}}{\partial \lambda}&=[f_1(\hat{x},x^{\star})]^{-1},\\
\frac{\partial x^{\star}}{\partial \lambda}&=[f_2(\hat{x},x^{\star})]^{-1},\\
\frac{\partial \hat{x}}{\partial \kappa} &=[f_3(\hat{x},x^{\star})]^{-1}\\
\frac{\partial x^{\star}}{\partial \kappa}&=[f_4(\hat{x},x^{\star})]^{-1}.
\end{align}
where the functions $f_1,f_2,f_3,f_4$ are given by:
\begin{align}
&\begin{aligned}
f_1(\hat{x},x^{\star})=
(\lambda+1)\left(\frac{1}{z}+\frac{1}{\hat{x}}\left[\frac{\hat{x}^{-l_1}+\hat{x}^{-l_2}}{l_1^{-1}z^{-l_1}+l_2^{-1}z^{-l_2}}-\frac{\hat{x}^{-l_1}-\hat{x}^{-l_2}}{l_1^{-1}z^{-l_1}-l_2^{-1}z^{-l_2}}\right]\right)
\label{derivative_1_principal_agent}\end{aligned}
\\&\begin{aligned}
f_2(\hat{x},x^{\star})=
(\lambda+1)\left(-\frac{1}{z}+\frac{1}{x^\star}\left[\frac{x^{\star-l_1}-x^{\star-l_2}}{l_1^{-1}z^{-l_1}-l_2^{-1}z^{-l_2}}+\frac{x^{\star-l_1}+x^{\star-l_2}}{l_1^{-1}z^{-l_1}+l_2^{-1}z^{-l_2}}\right]\right)
\end{aligned}
\\&\begin{aligned}
f_3(\hat{x},x^{\star})=
\frac{\kappa}{z}-\frac{1}{\hat{x}}\left(\kappa-z[l_1^{-1}+l_2^{-1}-1]\right)\left(\frac{\hat{x}^{-l_1}-\hat{x}^{-l_2}}{l_1^{-1}z^{-l_1}-l_2^{-1}z^{-l_2}}-\frac{\hat{x}^{-l_1}+\hat{x}^{-l_2}}{l_1^{-1}z^{-l_1}+l_2^{-1}z^{-l_2}}\right)&
\\-\frac{z}{\hat{x}}\frac{l_1z^{-l_1}-l_2z^{-l_2}}{l_1^{-1}z^{-l_1}+l_2^{-1}z^{-l_2}}\hspace{1 mm}&
\label{f3}
\end{aligned}
\\&\begin{aligned}
f_4(\hat{x},x^{\star})=-\frac{\kappa}{z}-\frac{1}{x^{\star}}\left(\kappa-z[l_1^{-1}+l_2^{-1}-1]\right)\left(\frac{x^{\star-l_1}-x^{\star-l_2}}{l_1^{-1}z^{-l_1}-l_2^{-1}z^{-l_2}}-\frac{x^{\star-l_1}+x^{\star-l_2}}{l_1^{-1}z^{-l_1}+l_2^{-1}z^{-l_2}}\right)&
\\+\frac{z}{x^{\star}}\frac{l_1z^{-l_1}-l_2z^{-l_2}}{l_1^{-1}z^{-l_1}+l_2^{-1}z^{-l_2}}\hspace{1 mm}&
\label{derivative_4_principal_agent}
\end{aligned}
\end{align}
\end{refproof}

\begin{refproof}[Proof of Corollary \ref{Corollary 3.4.}]

\noindent To prove Corollary \ref{Corollary 3.4.}, we firstly consider a control solution to the problem (\ref{external_objective_corr_3.4}) (which can by obtained using Theorem \ref{corollary_verification_theorem_zero_sum_degenerate}). Denote the optimal policy  
$u^{\star}_2\in\hspace{-2.0 mm} \underset{\hspace{3 mm }{u\in \mathcal{U}} }\arg\hspace{-1.4 mm}\sup J_2^{u}(t,x)$ where $u^{\star}_2=[\tau^{\star}_{2_j}, z^{\star}_{2_j}]_{j\geq 1}\in \mathcal{U}$ and the sets $\{\tau^{\star}_{2_j}\}_{j\in\mathbb{N}}$ and $\{z^{\star}_{2_j}\}_{j\in\mathbb{N}}$ are sequences of $\mathcal{F}_{\tau_j}-$measurable intervention times and intervention magnitudes respectively. Then by Remark \ref{Remark 4.3}, there exist constants $\hat{x}_2\in S$ and $x^{\star}_2\in S$ s.th $\hat{\tau}_{j+1}=\inf\{s>\tau_j;X^{\cdot,\hat{u} } (s)\geq x^{\star}_p \}\wedge \tau_S$ and $\hat{z}=\hat{x}_p-x^{\star}_p$. Hence, by setting $x^{\star}=x^{\star}_2$ and $\hat{x}=\hat{x}_2$ in Theorem \ref{Theorem 3.1.} we immediately deduce the result after applying the theorem.
\end{refproof}

Corollary \ref{Corollary 3.4.} demonstrates that the results can be applied to any pair of impulse control problems so that the cost parameters can be fixed so as to change the optimal control to match that of some other external objective function.

\section{Conclusion}

In this paper, we performed an analysis of the effects of imposing transaction cost on consumption behaviour. The inclusion of a transaction cost precludes agent behaviour policies in which the agent makes purchases continuously, hence in the model studied in this paper, the agent's behaviour is modelled using impulse control. In particular, we studied the effect of the transaction cost parameters on the consumption policy for an agent whose utility is given by a power utility function. The results of the paper provide a full characterisation of the parameters of transaction costs that sufficiently distorts the incentives of a rational agent so that the agent finds it optimal to adopt a consumption pattern that maximises the Principal's objective. Indeed, this paper describes for the first time, a Principal-Agent model with impulse control. Although the results of the paper are studied within the context of a liquidity-consumption problem, the results are broadly applicable. As described in Corollary \ref{Corollary 3.4.}, the results can be applied to any pair of impulse control problems so that the cost parameters can be fixed so as to change the optimal control to match that of some other external objective function.       \

An interesting avenue for future research is the effect of transaction costs on general impulse control problems in addition to specialised Principal objectives such as risk-minimisation and regime-dependant behaviour. 

\bibliographystyle{siamplain}

\end{document}